\documentclass{article}
\usepackage{amsmath,amssymb}
\newtheorem{theorem}{Theorem}
\newcommand{\bt}{\begin{theorem}}
\newcommand{\et}{\end{theorem}}
\newtheorem{lemma}{Lemma}
\newcommand{\bl}{\begin{lemma}}
\newcommand{\el}{\end{lemma}}
\newcommand{\pf}{{\bf Proof}.\ }
\newcommand{\be}{\begin{eqnarray}}
\newcommand{\ee}{\end{eqnarray}}
\newcommand{\beq}{\begin{equation}}
\newcommand{\eeq}{\end{equation}}
\newcommand{\benum}{\begin{enumerate}}
\newcommand{\eenum}{\end{enumerate}}
\newcommand{\balg}{\begin{align*}}
\newcommand{\ealg}{\end{align*}}
\newcommand{\ba}{\begin{array}}
\newcommand{\ea}{\end{array}}
\newcommand{\eop}{$\square$\vspace{.3cm}}

\begin{document}
\title{On the ubiquity of Sidon sets}
\author{Melvyn B. Nathanson\thanks{Supported
in part by grants from the NSA Mathematical Sciences Program
and the PSC-CUNY Research Award Program.
This paper was written while the author was a visitor
at the (alas, now defunct) AT\&T Bell Laboratories in Murray Hill, New Jersey, 
an excellent research institution that split into AT\&T Research Labs 
and Lucent Bell Labs, and provided another instance of a whole being greater 
than the sum of its parts.}\\
Department of Mathematics\\
Lehman College (CUNY)\\
Bronx, New York 10468\\
Email: nathansn@alpha.lehman.cuny.edu}

\maketitle

\begin{abstract}
A Sidon set is a set $A$ of integers such that no integer 
has two essentially distinct representations 
as the sum of two elements of $A$.
More generally, for every positive integer $g$, 
a $B_2[g]$-set is a set $A$ of integers such that no integer 
has more than $g$ essentially distinct representations 
as the sum of two elements of $A$.
It is proved that almost all small sumsets of $\{1,2,\ldots,n\}$ 
are $B_2[g]$-sets, in the sense that if $B_2[g](k,n)$
denotes the number of $B_2[g]$-sets of cardinality $k$ contained in the 
interval $\{1,2,\ldots,n\}$, then 
$\lim_{n\rightarrow\infty} B_2[g](k,n)/\binom{n}{k} = 1$
if $k = o\left(n^{g/(2g+2)}\right).$
\end{abstract}

\section{Sidon sets}
Let $A$ be a nonempty set of positive integers.  
The {\em sumset} $2A$ is the set of all integers of the form 
$a_1 + a_2,$ where $a_1,a_2 \in A$.
The set $A$ is called a {\em Sidon set} if every element of $2A$ 
has a unique representation as the sum of two elements of $A$, 
that is, if 
\[
a_1,a_2,a_1',a_2' \in A
\]
and
\[
a_1 + a_2 = a_1' + a_2',
\]
and if 
\[
a_1 \leq a_2 \quad \text{ and }\quad  a_1' \leq a_2',
\]
then 
\[
a_1 = a_1'  \quad \text{ and }\quad  a_2 = a_2'.
\]
More generally, for positive integers $h$ and $g$, 
the {\em $h$-fold sumset} $hA$ is the set of all 
sums of $h$ not necessarily distinct elements of $A$.
The {\em representation function} $r_{A,h}(m)$ counts the number 
of representations of $m$ in the form
\[
m = a_1 + a_2 + \cdots + a_h,
\]
where
\[
a_i\in A \qquad\text{for all $i = 1,2,\ldots,h$,}
\]
and
\[
a_1 \leq a_2 \leq \cdots \leq a_h.
\]
The set $A$ is called a $B_h[g]$-set if every element of $hA$ 
has at most $g$ representations as the sum of $h$ elements of $A$, 
that is, if 
\[
r_{A,h}(m) \leq g 
\]
for every integer $m$.
In particular, a $B_2[1]$-set is a Sidon set,
and $B_h[1]$-sets are usually denoted $B_h$-sets.

Let $h \geq 2$.  Let $A$ be a nonempty set of integers, and $a \in A.$
Then $r_h(m+a) \geq r_{h-1}(m).$  
Therefore, if $r_{A,h-1}(m) > g$ for some $m \in (h-1)A,$ 
then $r_{A,h}(m+a) > g$ for every $a \in A.$
It follows that if $A$ is a $B_h[g]$-set, 
then $A$ is also a $B_{h-1}[g]$-set.
In particular, every $B_h$- set is also a $B_{h-1}$-set.

Let $A$ be a subset of $\{1,2,\ldots,n\}$, and let $|A|$ denote the
cardinality of $A$.  Then $hA \subseteq \{h,h+1,\ldots,hn\}$.
If $|A| = k$, then there are exactly $\binom{k+h-1}{h}$
ordered $h$-tuples of the form $(a_1,\ldots,a_h)$,
where $a_i \in A$ for all $i = 1,\ldots,h$ and $a_1 \leq \cdots \leq a_h$.
If $A$ is a $B_h[g]$-set and $|A|= k,$ then
\[
\frac{k^h}{h!} < \binom{k+h-1}{h} = \sum_{m\in hA} r_{A,h}(m) \leq g|hA| < ghn,
\]
and so
\[
|A| = k < cn^{1/h}
\]
for $c = \left( h!gh \right)^{1/h}.$
It follows that if $A$ is a ``large'' subset of $\{1,2,\ldots,n\},$ 
then $A$ cannot be a $B_h[g]$-set.
In this paper we prove that almost all "small" subsets of $\{1,2,\ldots,n\}$
are $B_2[g]$-sets and almost all "small" subsets of $\{1,2,\ldots,n\}$
are $B_h$-sets.

{\em Notation.}  If $\{u_n\}_{n=1}^{\infty}$ and $\{v_n\}_{n=1}^{\infty}$ 
are sequences and $v_n > 0$ for all $n$,
we write $u_n = o(v_n)$ if $\lim_{n\rightarrow\infty} u_n/v_n = 0$,
and $u_n = O(v_n)$ or $u_n \ll v_n$ if $|u_n| \leq c v_n$ 
for some $c > 0$ and all $n \geq 1$.  The number $c$ in this inequality is 
called the {\em implied constant}.

\section{Random small $B_2[g]$-sets}
We require the following elementary lemma.

\bl          \label{sidon:lemma:binom}
If $n \geq 1$ and $0 \leq j \leq k \leq n$, then
\[
\frac{\binom{n-j}{k-j}}{\binom{n}{k}} \leq \left(\frac{k}{n}\right)^j.
\]
\el

\pf
We have
\[
\frac{\binom{n-j}{k-j}}{\binom{n}{k}} 
= \frac{(n-j)!k!}{(k-j)!n!} = \prod_{i=0}^{j-1}\frac{k-i}{n-i}
\leq \left(\frac{k}{n}\right)^j
\]
since
\[
\frac{k-i}{n-i} \leq \frac{k}{n}
\]
for $i = 0,1,\ldots, n-1$.
\eop

\bt  \label{sidon:theorem:lowerbound}
For any positive integers $g, k$, and $n$, 
let $B_2[g](k,n)$ denote the number of $B_2[g]$-sets $A$
contained in $\{1,\ldots,n\}$ with $|A| = k.$  
Then
\[
B_2[g](k,n) > \binom{n}{k}\left( 1 - \frac{4k^{2g+2}}{n^g}\right).
\]
\et

\pf
Let $A$ be a subset of $\{1,2,\ldots,n\}$ of cardinality $k$.
If $A$ is not a $B_2[g]$-set, then there is an integer $m \leq 2n$ 
such that $r_{A,2}(m) > g,$
that is, $m$ has at least $g+1$ representations as the sum of two elements of $A$.
This means that the set $A$ contains $g+1$ integers $a_1,\ldots, a_{g+1}$
such that 
\[
1 \leq a_1 < \cdots < a_{g+1} \leq \frac{m}{2},
\]
and $A$ also contains the $g+1$ integers $m - a_i$
for $i = 1,\ldots, g+1.$
If $a_{g+1} < m/2,$ then
\[
|\{a_i,m-a_i : i = 1,\ldots, g+1\}| = 2g-2.
\]
If $a_{g+1} = m/2,$ then
\[
|\{a_i,m-a_i : i = 1,\ldots, g+1\}| = 2g-1.
\]
Therefore, for each integer $m$, the number of sets $A \subseteq \{1,\ldots,n\}$
such that $|A| = k$ and $r_{A,2}(m) \geq g+1$ is at most
\[
\binom{\left[\frac{m-1}{2}\right]}{g+1}\binom{n-2g-2}{k-2g-2}
+ \binom{\left[\frac{m-1}{2}\right]}{g}\binom{n-2g-1}{k-2g-1},
\]
and so
\begin{align*}
\binom{n}{k} - B_2[g](k,n) \leq 
& \sum_{m\leq 2n} \binom{\left[\frac{m-1}{2}\right]}{g+1}\binom{n-2g-2}{k-2g-2}\\
& \quad + \sum_{m\leq 2n} \binom{\left[\frac{m-1}{2}\right]}{g}\binom{n-2g-1}{k-2g-1}.
\end{align*}
Observing that 
\[
\sum_{m\leq 2n} \binom{\left[\frac{m-1}{2}\right]}{g+1}
< \sum_{m\leq 2n} \left(\frac{m}{2}\right)^{g+1} \leq 2n^{g+2}
\]
and
\[
\sum_{m\leq 2n} \binom{\left[\frac{m-1}{2}\right]}{g}
< 2n^{g+1}
\]
and applying Lemma~\ref{sidon:lemma:binom}, we obtain
\begin{align*}
1 - \frac{B_2[g](k,n)}{\binom{n}{k}} 
& \leq \sum_{m\leq 2n} 
\left( 
\frac{\binom{\left[\frac{m-1}{2}\right]}{g+1}\binom{n-2g-2}{k-2g-2}}{\binom{n}{k}}
+ \frac{\binom{\left[\frac{m-1}{2}\right]}{g}\binom{n-2g-1}{k-2g-1}}{\binom{n}{k}}
\right)  \\
& \leq \sum_{m\leq 2n} 
\binom{\left[\frac{m-1}{2}\right]}{g+1}\left(\frac{k}{n}\right)^{2g+2}
+ \sum_{m\leq 2n} 
\binom{\left[\frac{m-1}{2}\right]}{g}\left(\frac{k}{n}\right)^{2g+1} \\
& < 2n^{g+2}\left(\frac{k}{n}\right)^{2g+2} + 2n^{g+1}\left(\frac{k}{n}\right)^{2g+1}\\
& \leq \frac{4k^{2g+2}}{n^g}.
\end{align*}
This completes the proof.
\eop

\bt           \label{sidon:theorem:Bhg}
Let $\{k_n\}_{n=1}^{\infty}$ be a sequence of positive integers 
such that $k_n \leq n$ for all $n$ and 
\[
k_n = o\left(n^{g/(2g+2)}\right).
\]
Then
\[
\lim_{n\rightarrow \infty} \frac{B_2[g](k_n,n)}{\binom{n}{k_n}} = 1.
\]
\et

\pf
This follows immediately from Theorem~\ref{sidon:theorem:lowerbound}.
\eop

\bt            \label{sidon:theorem:sidon}
Let $B_2(k,n)$ denote the number of Sidon sets of cardinality $k$ contained in 
$\{1,\ldots,n\}.$  If $k_n = o\left(n^{1/4}\right),$ then
\[
\lim_{n\rightarrow\infty} \frac{B_2(k_n,n)}{\binom{n}{k_n}} = 1.
\]
\et

\pf
This follows immediately from Theorem~\ref{sidon:theorem:Bhg} with $g = 1.$
\eop

\bt           \label{sidon:theorem:allsmall}
Let $\{k_n\}_{n=1}^{\infty}$ be a sequence of positive integers 
such that $k_n \leq n$ for all $n$ and 
\[
k_n = o\left(n^{g/(2g+3)}\right).
\]
Then
\[
\lim_{n\rightarrow \infty}
\frac{\sum_{k\leq k_n}  B_2[g](k,n)  }{\sum_{k\leq k_n}\binom{n}{k}} = 1.
\]
\et

\pf
It suffices to show that
\[
\lim_{n\rightarrow \infty}
\frac{\sum_{k\leq k_n} \left( \binom{n}{k} - B_2[g](k,n)\right) }{\sum_{k\leq k_n} \binom{n}{k}}
 = 0,
\]
where $f(k)$ is defined in the proof of Theorem~\ref{sidon:theorem:lowerbound}.
If $a_1,\ldots, a_{\ell}, b_1,\ldots, b_{\ell}$ are positive real numbers and 
$B = \max(b_1,\ldots, b_{\ell}),$ then
\[
\frac{a_1 + \cdots + a_{\ell}}{b_1 + \cdots + b_{\ell}}
\leq \frac{a_1 + \cdots + a_{\ell}}{B}  
\leq \frac{a_1}{b_1} + \cdots + \frac{a_{\ell}}{b_{\ell}}.  
\]
This implies that
\begin{align*}
\frac{\sum_{k\leq k_n} \left( \binom{n}{k} - B_2[g](k,n)\right)}{\sum_{k\leq k_n}\binom{n}{k}} 
& \leq \sum_{k\leq k_n} \frac{4k^{2g+2}}{n^g} \\
& \leq \frac{4k_n^{2g+3}}{n^g} \\
& = 4\left( \frac{k_n}{n^{g/(2g+3)}} \right)^{2g+3}\\
& = o(n),
\end{align*}
and the proof is complete.
\eop

We can restate our results in the language of probability.
Let $\Omega$ be the probability space consisting of the $\binom{n}{k}$  
subsets of $\{1,\ldots,n\}$ of cardinality $k$, 
where the probability of choosing $A \in \Omega$ is $1/\binom{n}{k}$.
If $P_{h,g}(k,n)$ denotes the probability that a random set $A \in \Omega$
is a $B_h[g]$-set, then Theorem~\ref{sidon:theorem:Bhg} states that
\[
\lim_{n\rightarrow\infty} P_{2,g}(k_n,n) = 1
\]
if $k_n = o\left(n^{g/(2g+2)}\right).$
 
Similarly, Theorem~\ref{sidon:theorem:allsmall} states that 
if $k_n = o\left(n^{g/(2g+3)}\right)$ and if 
$P_{h,g}(k_n,n)$ denotes the probability that a random set $A \subseteq \{1,\ldots,n\}$
of cardinality $|A| \leq k_n$ is a $B_h[g]$-set, then 
\[
\lim_{n\rightarrow\infty} P_{2,g}(k_n,n) = 1.
\]

\section{Random small $B_h$-sets}
A set $A$ is called a $B_h$-set if $r_{A,h}(m) = 1$ for all $m \in hA.$
Let $B_h(k,n)$ denote the number of $B_h$-sets of cardinality $k$
contained in $\{1,\ldots,n\}.$
Since every set of integers is a $B_1$-set, and every $B_h$ set
is a $B_{h-1}$-set, it follows that
\[
\binom{n}{k} = B_1(k,n) \geq \cdots \geq B_{h-1}(k,n) \geq B_h(k,n) \geq \cdots
\]

We shall prove that almost all ``small'' subsets of $\{1,\ldots,n\}$ are
$B_h$-sets.  The method is similar to that used to prove 
Theorem~\ref{sidon:theorem:Bhg}, but, for $h \geq 3,$ 
we have to consider the possible dependence between different 
representations of an integer $n$ as the sum of $h$ elements of $A$.
This means the following:  
Let $(a_1,\ldots,a_h)$ and $(a_1',\ldots,a_h')$
be $h$-tuples of elements of $A$ such that
\[
a_1 + \cdots + a_h = a_1' + \cdots + a_h',
\]
\[
a_1 \leq \cdots \leq a_h,
\]
\[
a_1' \leq \cdots \leq a_h',
\]
and
\[
\{a_1,\ldots,a_h\} \cap \{a_1',\ldots,a_h' \} \neq \emptyset.
\]  
If $h = 2,$ then $a_i = a_i'$ for $i = 1,2,$
but if $h \geq 3$, then it is not necessarily true that 
$a_i = a_i'$ for all $i = 1,\ldots, h.$
For example, in the case $h=3$ we have $1 + 3 + 4 = 1+2+5$ 
but $(1,3,4) \neq (1,2,5)$.  In the case $h=5$ we have
$1+1+2+3+3 = 1+2+2+2+3$ but $(1,1,2,3,3) \neq (1,2,2,2,3),$
even though $\{1,1,2,3,3\} = \{1,2,2,2,3\}.$

Because of the lack of independence, we need a careful description
of a representation of $m$ as the sum of $h$ not necessarily 
distinct integers.  
We introduce the following notation.
Let $A$ be a set of positive integers.
Corresponding to each representation of $m$ in the form
\beq                  \label{sidon:reph}
m = a_1 + \cdots + a_h,
\eeq
where $a_1, \ldots, a_h \in A$ and $a_1 \leq \cdots \leq a_h$,
there is a unique triple 
\beq                                \label{sidon:triple}
(r,(h_j),(a_j')),
\eeq
where
\benum
\item[(i)]
$r$ is the number of distinct summands in this representation, 
\item[(ii)]
$(h_j) = (h_1,\ldots,h_r)$ is an ordered partition of $h$ into $r$ positive parts,
that is, an $r$-tuple of positive integers such that
\[
h = h_1 + \cdots + h_r,
\]
\item[(iii)]
$(a_j') = (a'_1,\ldots,a'_r)$ is an $r$-tuple of pairwise distinct 
elements of $A$ such that
\[
1 \leq a'_1 < \cdots < a'_r \leq m
\]
and
\[
\{a_1,\ldots,a_h\} = \{a'_1,\ldots,a'_r\},
\]
\item[(iv)]
\[
m = h_1a'_1 + \cdots + h_ra'_r,
\]
where each integer $a'_j$ occurs exactly $h_j$ times 
in the representation~(\ref{sidon:reph}).
\eenum
There is a one-to-one correspondence between distinct
representations of an integer $m$ in the form~(\ref{sidon:reph})
and triples of the form~(\ref{sidon:triple}).
Moreover, for each $r$ and $m$, the integer $a'_r$ is completely 
determined by the ordered partition $(h_j)$ of $h$ and the
$(r-1)$-tuple $(a'_1,\ldots,a'_{r-1})$.  Therefore, for positive
integers $m$ and $r$, the number of triples 
of the form~(\ref{sidon:triple}) does not exceed $\pi_r(h)m^{r-1},$
where $\pi_r(h)$ is number of ordered partitions of $h$ into exactly $r$ positive parts.

\bt       \label{sidon:theorem:Bhdiff}
Let $h \geq 2.$  For all positive integers $k \leq n,$
\[
B_{h-1}(k,n) - B_h(k,n) \ll \binom{n}{k} \frac{k^{2h}}{n}
\]
and
\[
B_h(k,n) \gg \binom{n}{k} \left( 1 - \frac{k^{2h}}{n} \right). 
\]
where the implied constants depend only on $h$.
\et

\pf
Let $A$ be a $B_{h-1}$-set contained in $\{1,\ldots,n\}.$ 
Then $hA \subseteq \{h,h+1,\ldots,hn\}.$
If $m \in hA$ and $m$ has two distinct representations
as the sum of $h$ elements of $A$, then there exist positive integers 
$r_1$ and $r_2$ and triples
\beq           \label{sidon:2triples}
(r_1,(h_{1,j}),(a'_{1,j}))  \qquad\text{and}\qquad (r_2,(h_{2,j}),(a'_{2,j}))
\eeq
such that, for $i = 1$ and 2, we have
\[
m = \sum_{j=1}^{r_i} h_{i,j}a'_{i,j},
\]
\[
h = \sum_{j=1}^{r_i} h_{i,j},
\]
and
\[
1 \leq a'_{i,1} < \cdots < a'_{i,r_i} \leq m.
\]
The number of pairs of triples of the form~(\ref{sidon:2triples})
for fixed positive integers $m$, $r_1 \leq h$, and $r_2 \leq h$ 
is at most
\[
\pi_{r_1}(h)m^{r_1-1}\pi_{r_2}(h)m^{r_2-1} \ll m^{r_1+r_2-2},
\]
where the implied constant depends only on $h$.
Moreover, since $A$ is a $B_{h-1}$-set, 
no number can have two representations 
as the sum of $h-1$ elements of $A$.  This implies that
\[
\{a'_{1,1}, a'_{1,2},\ldots,a'_{1,r_1}\}
\cap \{a'_{2,1}, a'_{2,2},\ldots,a'_{2,r_2}\}
= \emptyset,
\]
and so the set
\[
\{a'_{i,j} : i = 1,2 \text{ and } j = 1,\ldots,r_i\}
\]
contains exactly $r_1+r_2$ elements of $A$.
Therefore, given positive integers $r_1 \leq h$, $r_2 \leq h$
and $m \leq hn$, there are
\[
\ll m^{r_1+r_2-2}\binom{n-r_1-r_2}{k-r_1-r_2}
\]
sets $A$ for which $m$ has two representations
as the sum of $h$ elements of $A$, and in which one representation
uses $r_1$ distinct integers and the other representation uses $r_2$
distinct integers.  Summing over $m \leq hn,$ we obtain
\[
\ll n^{r_1+r_2-1}\binom{n-r_1-r_2}{k-r_1-r_2}.
\]
Applying Lemma~\ref{sidon:lemma:binom}, we obtain
\begin{align*}
B_{h-1}(k,n) - B_{h}(k,n)
& \ll \binom{n}{k}
\sum_{r_1,r_2 \leq h} 
\frac{n^{r_1+r_2-1}\binom{n-r_1-r_2}{k-r_1-r_2}}{\binom{n}{k}} \\
& \ll \binom{n}{k}\sum_{r_1,r_2 \leq h} n^{r_1+r_2-1}\left(\frac{k}{n}\right)^{r_1+r_2}  \\
& = \binom{n}{k}\sum_{r_1,r_2 \leq h} \frac{k^{r_1+r_2}}{n}  \\
& \ll \binom{n}{k}\frac{k^{2h}}{n}.
\end{align*}
It follows that
\begin{align*}
\binom{n}{k} - B_h(k,n)
& = B_1(k,n) - B_h(k,n)  \\
& = \sum_{j=2}^h \left(B_{j-1}(k,n) - B_j(k,n)\right)  \\
& \ll \sum_{j=2}^h \binom{n}{k}\frac{k^{2j}}{n} \\
& \ll \binom{n}{k}\frac{k^{2h}}{n},
\end{align*}
and so
\[
B_h(k,n) \gg \binom{n}{k} \left( 1 - \frac{k^{2h}}{n} \right). 
\]
This completes the proof.
\eop

\bt       \label{sidon:theorem:Bhsets}
Let $B_h(k,n)$ denote the number of $B_h$-sets $A$ contained in $\{1,\ldots,n\}$
with $|A| = k.$  
Let $\{k_n\}_{n=1}^{\infty}$ be a sequence of positive integers 
such that 
\[
k_n = o\left(n^{1/2h}\right).
\]
Then
\[
\lim_{n\rightarrow\infty} \frac{B_h(k_n,n)}{\binom{n}{k_n}} = 1.
\]
\et

\pf
By Theorem~\ref{sidon:theorem:Bhdiff},
\[
1 \geq \frac{B_h(k_n,n)}{\binom{n}{k_n}} \gg 1 - \left( \frac{k_n}{n^{1/2h}} \right)^{2h},
\]
and so
\[
\lim_{n\rightarrow\infty} \frac{B_h(k_n,n)}{\binom{n}{k_n}} = 1.
\]
This completes the proof.
\eop

\section{Remarks added in proof}
A variant of Theorem~\ref{sidon:theorem:sidon} appears in 
Nathanson~\cite[p. 37, Exercise 14]{nath96bb}.
Godbole, Janson, Locantore, and Rapoport~\cite{godb-jans-loca-rapo99} used
probabilistic methods to obtain a converse of Theorem~\ref{sidon:theorem:Bhsets}.
They proved that if
\[
\lim_{n\rightarrow\infty} \frac{k_n}{n^{1/2h}} = \infty,
\]
then
\[
\lim_{n\rightarrow\infty} \frac{B_h(k_n,n)}{\binom{n}{k_n}} = 0.
\]
They also analyzed the threshold behavior of $B_h(k,n)$, and proved that if
\[
\lim_{n\rightarrow\infty} \frac{k_n}{n^{1/2h}} = \Lambda > 0,
\]
then 
\[
\lim_{n\rightarrow\infty} \frac{B_h(k_n,n)}{\binom{n}{k_n}} = e^{-\lambda},
\]
where $\lambda = \kappa_h\Lambda^{2h}.$

It is natural to conjecture that analogous results 
hold for the function $B_h[g](k,n)$, namely, if 
\[
\lim_{n\rightarrow\infty} \frac{k_n}{n^{g/(gh+h)}} = 0,
\]
then
\[
\lim_{n\rightarrow\infty} \frac{B_h(k_n,n)}{\binom{n}{k_n}} = 1,
\]
and if 
\[
\lim_{n\rightarrow\infty} \frac{k_n}{n^{g/(gh+h)}} = \infty,
\]
then
\[
\lim_{n\rightarrow\infty} \frac{B_h(k_n,n)}{\binom{n}{k_n}} = 0.
\]
It should also be possible to describe the threshold behavior 
of $B_h(k_n,n)/\binom{n}{k_n}$ in the case 
\[
\lim_{n\rightarrow\infty} \frac{B_h(k_n,n)}{\binom{n}{k_n}} = \Lambda > 0.
\]

\providecommand{\bysame}{\leavevmode\hbox to3em{\hrulefill}\thinspace}
\providecommand{\MR}{\relax\ifhmode\unskip\space\fi MR }
\providecommand{\MRhref}[2]{%
  \href{http://www.ams.org/mathscinet-getitem?mr=#1}{#2}
}
\providecommand{\href}[2]{#2}


\begin{thebibliography}{1}

\bibitem{godb-jans-loca-rapo99}
A.~P. Godbole, S.~Janson, {N. W. Locantore, Jr.}, and R.~Rapoport, \emph{Random
  {S}idon sequences}, J. Number Theory \textbf{75} (1999), 7--22.

\bibitem{nath96bb}
M.~B. Nathanson, \emph{Additive number theory: Inverse problems and the
  geometry of sumsets}, Graduate Texts in Mathematics, vol. 165,
  Springer-Verlag, New York, 1996.

\end{thebibliography}
\end{document}